\documentclass{article}

\PassOptionsToPackage{square, numbers, sort&compress}{natbib}

\usepackage[final]{neurips_2023_ml4ps}

\usepackage[utf8]{inputenc} 
\usepackage[T1]{fontenc}    
\usepackage{hyperref}       
\usepackage{url}            
\usepackage{booktabs}       
\usepackage{amsmath, amsfonts}       
\usepackage{nicefrac}       
\usepackage{microtype}      
\usepackage{xcolor}         
\usepackage{multirow} 
\usepackage{graphicx}
\usepackage{caption, subcaption}

\newcommand{\eref}[1]{\mbox{\rm(\ref{#1})}}

\newcommand{\set}[1]{\left\{ #1 \right\} }
\newcommand{\norm}[1]{\left\| #1 \right\|}

\newcommand{\real}{\mathbb{R}}

\newcommand {\cD} {{\cal D}}

\newcommand {\cL} {{\cal L}}

\newcommand{\bzero}{\boldsymbol{0}}

\newcommand{\bg}{\boldsymbol{g}}
\newcommand{\bh}{\boldsymbol{h}}

\newcommand{\br}{\boldsymbol{r}}

\newcommand{\bx}{\boldsymbol{x}}

\newcommand{\bA}{\boldsymbol{A}}
\newcommand{\bB}{\boldsymbol{B}}
\newcommand{\bC}{\boldsymbol{C}}

\newcommand{\bI}{\boldsymbol{I}}

\newcommand{\bX}{\boldsymbol{X}}

\newcommand{\hbx}{\widehat{\boldsymbol{x}}}

\newcommand{\bmu}{{\boldsymbol{\mu}}}

\title{Nonlinear-manifold reduced order models with domain decomposition}

\author{%
  Alejandro N.~Diaz \\
  Rice University\\
  Houston, TX 77005 \\
  \texttt{and5@rice.edu} \\
  \And
  Youngsoo Choi \\
  Lawrence Livermore National Laboratory \\
  Livermore, CA 94550 \\
  \texttt{choi15@llnl.gov} \\
  \And
  Matthias Heinkenschloss \\
  Rice University\\
  Houston, TX 77005 \\
  \texttt{heinken@rice.edu} \\
}

\begin{document}

\maketitle

\begin{abstract}
A nonlinear-manifold reduced order model (NM-ROM) is a great way of incorporating underlying physics principles into a neural network-based data-driven approach. We combine NM-ROMs with domain decomposition (DD) for efficient computation. NM-ROMs offer benefits over linear-subspace ROMs (LS-ROMs) but can be costly to train due to parameter scaling with the full-order model (FOM) size. To address this, we employ DD on the FOM, compute subdomain NM-ROMs, and then merge them into a global NM-ROM. This approach has multiple advantages: parallel training of subdomain NM-ROMs, fewer parameters than global NM-ROMs, and adaptability to subdomain-specific FOM features. Each subdomain NM-ROM uses a shallow, sparse autoencoder, enabling hyper-reduction (HR) for improved computational speed. In this paper, we detail an algebraic DD formulation for the FOM, train HR-equipped NM-ROMs for subdomains, and numerically compare them to DD LS-ROMs with HR. Results show a significant accuracy boost, on the order of magnitude, for the proposed DD NM-ROMs over DD LS-ROMs in solving the 2D steady-state Burgers' equation.
\end{abstract}

\section{Introduction}\label{sec:introduction}
In science and engineering, complex tasks often involve repeatedly simulating a large-scale, parameterized, nonlinear system referred to as the full-order model (FOM). Ensuring high fidelity requires a high-dimensional model, leading to significant computational costs and lengthy simulations. As a result, tasks like design optimization become impractical for large-scale problems. Model reduction offers a solution by replacing the FOM with a computationally efficient, low-dimensional model called a reduced-order model (ROM). This ROM approximates the FOM's behavior with adjustable accuracy, making it suitable for many-query applications. However, construction of accurate and computationally efficient ROMs poses challenges.
To address them, we integrate the nonlinear-manifold ROM (NM-ROM) approach with an algebraic domain-decomposition (DD) framework.

Various model reduction methods have been integrated with DD, like reduced basis elements (RBE) \cite{YMaday_EMRonquist_2002a,
YMaday_EMRonquist_2004a,
LIapichino_AQuarteroni_GRozza_2012a,
PFAntonietti_PPacciarini_AQuarteroni_2016a,
JLEftang_DBPHuynh_DJKnezevic_EMRonquist_ATPatera_2012a,
DBPHuynh_DJKnezevic_ATPatera_2013a,JLEftang_ATPatera_2013a,
LIapichino_AQuarteroni_GRozza_2012a}, and the alternating Schwarz method \cite{MBuffoni_HTelib_AIollo_2009a,
JBarnett_ITezaur_AMota_2022a,
KSmetana_TTaddei_2022a,AIollo_GSambataro_TTaddei_2022a}. 
However, they are often specialized to specific problems, dealing with the physical domain at the PDE level. 
In contrast, the authors in \cite{CHoang_YChoi_KCarlberg_2021a} take an algebraic approach by decomposing the FOM at the discrete level and computing linear-subspace ROMs (LS-ROMs) for each subdomain. While LS-ROMs work well in many cases \cite{BHaasdonk_2017a,
AQuarteroni_AManzoni_FNegri_2016a, 
MHinze_SVolkwein_2005a,
MGubisch_SVolkwein_2017a,
cheung2023local,
copeland2022reduced,carlberg2018conservative,
ACAntoulas_2005a,PBenner_TBreiten_2017a, 
ACAntoulas_CABeattie_SGugercin_2020a,CGu_2011a,PBenner_TBreiten_2015a,
AJMayo_ACAntoulas_2007a,
ACAntoulas_IVGosea_ACIonita_2016a,
IVGosea_ACAntoulas_2018a,
choi2021space,
kim2021efficient,
choi2019space}, 
it is well known that advection-dominated problems and problems with sharp gradients cannot be well-approximated using low-dimensional linear subspaces. 
These problems are said to have slowly decaying Kolmogorov $n$-width \cite{MOhlberger_SRave_2016a}.
Recent approaches, such as nonlinear-manifold ROMs (NM-ROMs), address these problems by nonlinearly approximating the FOM in a low-dimensional nonlinear manifold. 
This is typically achieved through training an autoencoder on FOM snapshot data (e.g., \cite{
KKashima_2016a,
DHartman_LKMestha_2017a,
KLee_KTCarlberg_2020a,
YKim_YChoi_DWidemann_TZohdi_2022a,
kim2020efficient}).
However, training of NM-ROMs is expensive.
Indeed, in the monolithic single-domain case, the high-dimensionality of the FOM training data results in a large number of neural network (NN) parameters requiring training. 
In \cite{JBarnett_CFarhat_YMaday_2023a} this cost issue was mitigated by first computing a  low dimensional proper orthogonal decomposition
(POD) model, and then using a NN to train the coefficients in this POD.
Instead, we integrate an autoencoder framework with DD.
By coupling NM-ROM with DD, one can compute FOM training data on subdomains, thus reducing the dimensionality of subdomain NM-ROM training data, resulting in fewer parameters that need to be trained per subdomain NM-ROM. 

We also note that couplings of NNs and DD for solutions of partial differential equations (PDEs) have been considered in previous work (e.g., \cite{KLi_KTang_TWu_QLiao_2020a, 
WLi_XXiang_YXu_2020a, 
QSun_XXu_HYi_2023a, 
SLi_YXia_YLiu_QLiao_2023a}).
However, these approaches use deep learning to solve a PDE by representing its solution as a NN and minimizing a corresponding physics-informed loss function.
In contrast, our work uses autoencoders to reduce the dimensionality of an existing numerical model. 
The autoencoders are pretrained in an {\it offline} stage to find low-dimensional representations of FOM snapshot data, and used in an {\it online} stage to significantly reduce the computational cost and runtime of numerical simulations. 
Our work is the first to couple autoencoders with DD in the reduced-order modeling context.

Here, we extend the work of \cite{CHoang_YChoi_KCarlberg_2021a} on DD LS-ROM and integrate NM-ROM with hyper-reduction (HR) using shallow, sparse autoencoders discussed in \cite{YKim_YChoi_DWidemann_TZohdi_2022a}. 
We incorporate the NM-ROM approach into this framework because of its success when applied to problems with slowly decaying Kolmogorov $n$-width. 
DD allows one to compute FOM training snapshots on subdomains, thus reducing the dimensionality of subdomain NM-ROM training data, resulting in fewer parameters that need to be trained per subdomain NM-ROM.
We use wide, shallow, and sparse autoencoder architecture, which allows HR to be efficiently applied, thus reducing the complexity caused by nonlinearity and yielding computational speedup. 
Additionally, we modify the wide, shallow, and sparse architecture used in \cite{YKim_YChoi_DWidemann_TZohdi_2022a} to also include a sparsity mask for the encoder input layer as well as the decoder output layer. 
The proposed DD NM-ROM approach is compared with DD LS-ROM on the 2D Burgers' equation.
\section{DD full order model}\label{sec:dd_fom}
First consider the monolithic, single-domain FOM written as a residual equation
\begin{equation}\label{eq:fom_monolithic}
    \br(\bx; \bmu) = \bzero,
\end{equation}
where 
$\bx\in \real^{N_x}$ is the state, 
$\bmu \in \cD \subset \real^{N_\mu}$ is a parameter, and
$\br:\real^{N_x}\times \real^{N_\mu}\to \real^{N_x}$ is the residual function. 
FOMs of the form \eref{eq:fom_monolithic} typically arise from discretizations of partial differential equations (PDEs).
One can reformulate \eref{eq:fom_monolithic} into a DD formulation by partitioning the residual equation into $n_\Omega$ systems of equations (so-called {\it algebraic} subdomains), coupling them via {\it compatibility constraints}, and converting the systems of equations into a least-squares problem, resulting in
\begin{equation}\label{eq:dd_fom}
\min_{(\bx_i^\Omega, \bx_i^\Gamma), i=1, \dots, n_\Omega} \quad 
\frac{1}{2}\sum_{i=1}^{n_\Omega} \norm{\br_i\left(\bx_i^\Omega, \bx_i^\Gamma; \bmu \right)}_2^2, \quad
{\rm s.t.} \quad \sum_{i=1}^{n_\Omega} \bA_i \bx_i^\Gamma = \bzero,
\end{equation} where
$\bx_i^\Omega\in \real^{N_i^\Omega}$,
$\bx_i^\Gamma\in \real^{N_i^\Gamma}$,
$\br_i:\real^{N_i^\Omega}\times \real^{N_i^\Gamma}\times \cD \to \real^{N_i^r}$, and
$\bA_i\in \set{-1, 0, 1}^{N_a\times N_i^\Gamma}$
are the $i$-th subdomain interior-state, interface-state, 
residual function, and compatibility constraint matrix, respectively.
The sparsity pattern of the monolithic residual function $\br$ determines the structure of the subdomain residual functions $\br_i$, as well as the decomposition of the state $\bx$ into subdomain states $(\bx_i^\Omega, \bx_i^\Gamma)$. 
The interior-states $\bx_i^\Omega$ are those that are \emph{only} used to compute the residual $\br_i$ in the $i$-th subdomain, whereas the interface-states $\bx_i^\Gamma$ are also used in the residual computation of neighboring subdomains. 
The equality constraint determined by $\bA_i$ enforces equality on the overlapping interface states. For further details, see \cite[Sec. 2]{ANDiaz_YChoi_MHeinkenschloss_2023a} or \cite[Sec. 2]{CHoang_YChoi_KCarlberg_2021a}.
\section{DD nonlinear-manifold reduced order model}\label{sec:dd_nmrom}
For each subdomain $i \in \set{1, \dots, n_\Omega}$, 
let $\bg_i^\Omega: \real^{n_i^\Omega} \to \real^{N_i^\Omega}$, $n_i^\Omega \ll N_i^\Omega$, and 
$\bg_i^\Gamma: \real^{n_i^\Gamma} \to \real^{N_i^\Gamma}$, $n_i^\Gamma \ll N_i^\Gamma$,
be decoders such that 
$\bx_i^\Omega \approx \bg_i^\Omega(\hbx_i^\Omega)$ and
$\bx_i^\Gamma \approx \bg_i^\Gamma(\hbx_i^\Gamma)$.
Also let $\bB_i \in \set{0,1}^{N_i^B \times N_i^r}$, $N_i^B \leq N_i^r$, denote a row-sampling matrix for collocation HR, and let
$\bC\in \real^{n_C\times N_{\overline{A}}}$, $\; n_C \ll N_a$, be a Gaussian test matrix.
The DD NM-ROM is evaluated by solving
\begin{equation}\label{eq:dd_nmrom}
\min_{(\hbx_i^\Omega, \hbx_i^\Gamma), i=1, \dots, n_\Omega} \quad 
\frac{1}{2}\sum_{i=1}^{n_\Omega} \norm{\bB_i \br_i\left(\bg_i^\Omega\left(\hbx_i^\Omega\right), \bg_i^\Gamma\left(\hbx_i^\Gamma \right)\right)}_2^2, \quad
{\rm s.t.} \quad  \sum_{i=1}^{n_\Omega} \bC\bA_i \bg_i^\Gamma(\hbx_i^\Gamma) = \bzero.
\end{equation}

If HR is not applied (i.e., $\bB_i=\bI$ in \eref{eq:dd_nmrom}), the ROM's computational savings are limited because evaluation of residuals
$\big(\hbx_i^\Omega, \hbx_i^\Gamma \big) \to$  
$\big( \bg_i^\Omega(\hbx_i^\Omega), \bg_i^\Gamma(\hbx_i^\Gamma) \big) \to$ 
$\br_i\big(\bg_i^\Omega\big(\hbx_i^\Omega\big), \bg_i^\Gamma\big(\hbx_i^\Gamma \big) \big)$ 
scales with the size $N_i^\Omega$ and $N_i^\Gamma$ of the FOM. 
Thus, HR is applied to decrease the computational complexity caused by the nonlinearity of $\br_i$, and increase the computational speedup. 
We use \cite[Algo. 3]{KTCarlberg_CFarhat_JCortial_DAmsallem_2013a} to greedily compute a row sampling matrix $\bB_i$ for collocation HR. 
The application of HR to the decoders $\bg_i^\Omega$ and $\bg_i^\Gamma$ is discussed further in Sec. \ref{sec:architecture}.
Following \cite{CHoang_YChoi_KCarlberg_2021a}, we apply a Gaussian test matrix $\bC\in \real^{n_C\times N_a}$, $\; n_C \ll N_{a}$, to convert the compatibility constraints into a
so-called ``weak compatibility constraint", which decreases the number of constraints to avoid making the DD ROM over-determined.

The DD FOM \eref{eq:dd_fom} and DD NM-ROM \eref{eq:dd_nmrom} are solved using an inexact Lagrange-Newton sequential quadratic programming (SQP) solver, where the Hessian of the Lagrangian is replaced with a Gauss-Newton approximation.
This avoids computation of second order derivatives of residuals and constraints in \eref{eq:dd_nmrom},
but still achieves good convergence for \eref{eq:dd_fom} and \eref{eq:dd_nmrom}. 
For further details, see \cite{ANDiaz_YChoi_MHeinkenschloss_2023a}.

The DD NM-ROM \eref{eq:dd_nmrom} formulation has several benefits. Training, i.e., computation of the $\bg_i^\Omega$ and
$\bg_i^\Gamma$ is local, involves few parameters, and can be done in parallel.
The ROMs can be adjusted to localized features of the problem, which may result in smaller ROMs. 
Parallelization can be used to speed up ROM computation/training and ROM execution.
\subsection{NM-ROM architecture and training}\label{sec:architecture}

We use single-layer, wide, and sparse decoders with smooth activation functions to represent the maps $\bg_i^\Omega$ and $\bg_i^\Gamma$. 
The corresponding encoders, denoted
$\bh_i^\Omega$ and
$\bh_i^\Gamma$,
are also single-layer, wide, and sparse.
Shallow networks are used for computational efficiency; fewer layers correspond to fewer repeated matrix-vector multiplications when evaluating the decoders. 
The shallow depth necessitates a wide network to maintain enough expressiveness for use in NM-ROM. 
Smooth activations (i.e., swish) are used to ensure that $\bg_i^\Omega$ and $\bg_i^\Gamma$ are continuously differentiable. 
Normalization and de-normalization layers are also applied at the encoder input and decoder output layers, respectively.

Sparsity is applied at the decoder output layer so that HR can be applied. 
The sparsity allows one to compute a \emph{subnet}, which only keeps track of the hidden nodes required to compute the output nodes that remain after HR. 
Further details can be found in \cite[Sec. 3.2]{YKim_YChoi_DWidemann_TZohdi_2022a}, 
\cite[Sec. 5.3]{ANDiaz_YChoi_MHeinkenschloss_2023a}.
We also apply a sparsity mask to the encoder input layer so that the autoencoders are symmetric across the latent layer. The sparsity pattern has a tri-banded structure inspired by 2D finite difference stencils, where the number of nonzeros per band and the separation between bands are hyper-parameters. 

To train the autoencoders, we first generate FOM snapshots in an {\it offline} stage by solving \eref{eq:dd_fom} at parameters $\set{\bmu_\ell}_{\ell=1}^M$, and collect interior- and interface-state snapshot datasets $\bX_i^\Omega \in \real^{N_i^\Omega\times M}$ and $\bX_i^\Gamma\in \real^{N_i^\Gamma\times M}.$
Alternatively, one can solve the monolithic FOM \eref{eq:fom_monolithic} at each $\bmu_\ell$ and restrict the corresponding states $\bx(\bmu_\ell)$ to interior-states $\bx_i^\Omega(\bmu_\ell)$ and interface-states $\bx_i^\Gamma(\bmu_\ell)$ for each subdomain.
We use the latter approach.
The autoencoders $(\bh_i^\Omega, \bg_i^\Omega)$ and $(\bh_i^\Gamma, \bg_i^\Gamma)$ are then trained in parallel by minimizing the respective MSE losses
\begin{equation}\label{eq:MSE_loss}
\hspace*{-2ex}
    \cL_i^\Omega = \frac{1}{M}\sum_{\ell=1}^{M} \norm{\bx_i^{\Omega}(\bmu_\ell)-\bg_i^\Omega(\bh_i^\Omega(\bx_i^{\Omega}(\bmu_\ell)))}_2^2, \;
    \cL_i^\Gamma = \frac{1}{M}\sum_{\ell=1}^{M} \norm{\bx_i^{\Gamma}(\bmu_\ell)-\bg_i^\Gamma(\bh_i^\Gamma(\bx_i^{\Gamma}(\bmu_\ell)))}_2^2
\end{equation}
for each subdomain $i = 1, \ldots, n_\Omega$.
The snapshots undergo a random 90-10 split for training and validation, and the MSE loss is minimized using the Adam optimizer over $2000$ epochs with a batch size of $32$. 
We also apply early stopping \cite{LPrechelt_1998a} with a patience of $300$ and reduce the learning rate on plateau with an initial learning rate of $10^{-3}$. 
The implementation was done in PyTorch and used the PyTorch Sparse and SparseLinear packages.
\section{Numerical experiment: 2D Burgers' equation}\label{sec:numerics}
We compare the DD LS-ROM of \cite{CHoang_YChoi_KCarlberg_2021a} and the proposed DD NM-ROM with HR for the 2D steady-state Burgers equation. The DD LS-ROM can be regarded as a specific instance within the realm of DD NM-ROMs, where the encoders and decoders defined in Equation (4) are exchanged for linear operators derived through singular value decomposition. We compute the relative error as
\begin{equation}\label{eq:relative_error}
    e = \left(\frac{1}{n_\Omega}\sum_{i=1}^{n_\Omega} \Big( \norm{\bx_i^\Omega-\bg_i^\Omega(\hbx_i^\Omega)}_2^2+\norm{\bx_i^\Gamma-\bg_i^\Gamma(\hbx_i^\Gamma)}_2^2\Big) / \Big(\norm{\bx_i^\Omega}_2^2+\norm{\bx_i^\Gamma}_2^2\Big)\right)^{1/2}.
\end{equation}
All training and computations were performed on the Lassen machine at Lawrence Livermore National Laboratory, which consists of an IBM Power9 processor with NVIDIA V100 (Volta) GPUs, clock speed between 2.3-3.8 GHz, and 256 GB DDR4 memory. 
The code can be found at https://anonymous.4open.science/r/DDNMROM\_NeurIPS-4160/.

The implementation was done sequentially, but to highlight potential advantages of a parallel implementation, the reported wall clock time for computing subdomain-specific quantities for the SQP solver is taken to be the largest wall clock time incurred among all subdomains. 
The wall clock time for the remaining steps of the SQP solver is set to the overall wall clock time.

\begin{figure}[ht]
     \centering
     \begin{subfigure}[b]{0.475\textwidth}
         \centering
         \includegraphics[width=0.48\textwidth]{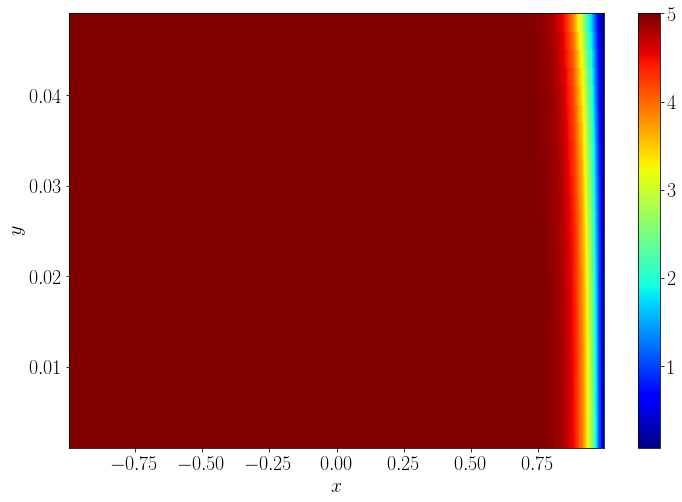} \hfill
         \includegraphics[width=0.48\textwidth]{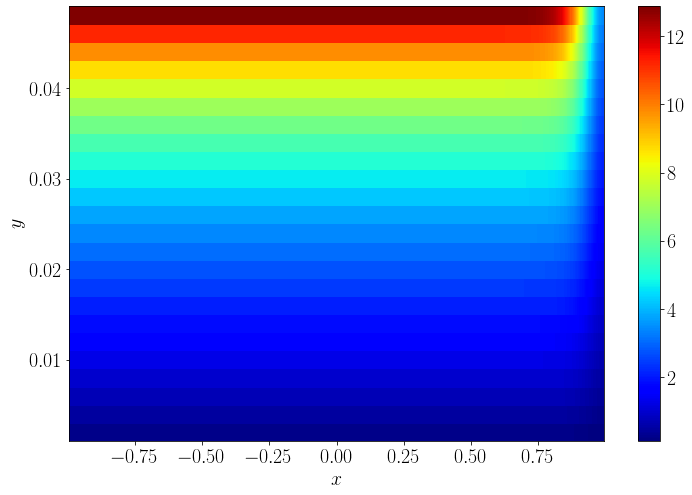}
         \caption{$(a, \lambda) = (1, 25)$.}
         \label{fig:uv_1}
     \end{subfigure}
     \hfill
     \begin{subfigure}[b]{0.475\textwidth}
         \centering
         \includegraphics[width=0.48\textwidth]{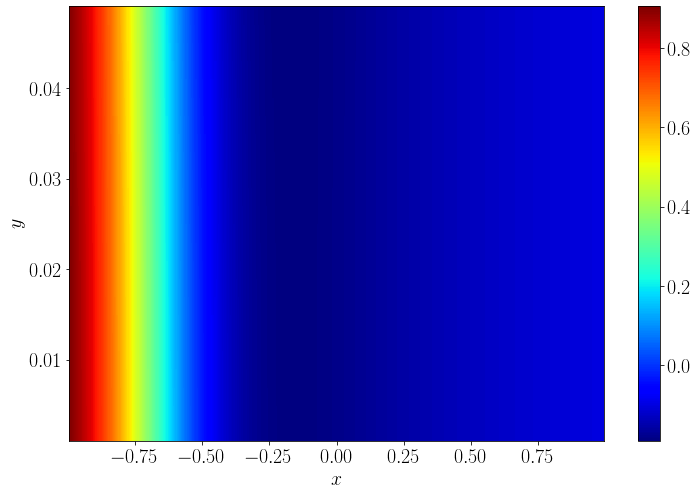} \hfill
         \includegraphics[width=0.48\textwidth]{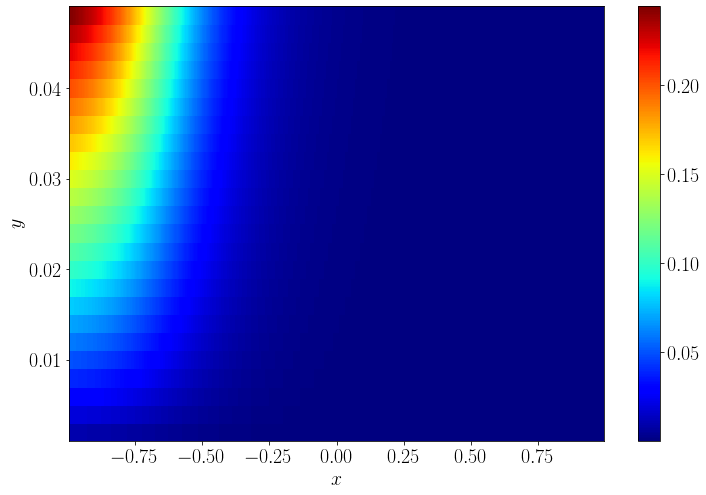}
         \caption{$(a, \lambda) = (10^4, 5)$.}
         \label{fig:uv_2}
     \end{subfigure}
        \caption{FOM $u$ and $v$ components for different $(a, \lambda)$. The distance of the shock from the left boundary and its steepness are determined by $a$ and $\lambda$, respectively.}
        \label{fig:uv_states}
\end{figure}

We consider the 2D steady-state Burgers' equation 
\begin{align}\label{eq:burgers_pde_2d}
    u \frac{\partial u}{\partial x} + v \frac{\partial u}{\partial y} &= \nu \left(\frac{\partial^2 u}{\partial x^2} + \frac{\partial^2 u}{\partial y^2}\right),
    &
    u \frac{\partial v}{\partial x} + v \frac{\partial v}{\partial y} &= \nu \left(\frac{\partial^2 v}{\partial x^2} + \frac{\partial^2 v}{\partial y^2}\right)
\end{align}
for $(x, y)\in [-1, 1] \times [0, 0.05]$ with viscosity $\nu=0.1$.
As in \cite{CHoang_YChoi_KCarlberg_2021a}, we use the exact solution
$u_{ex} = - 2\nu  \frac{\partial}{\partial x} \psi\, / \psi$,\;
$v_{ex} = - 2\nu  \frac{\partial}{\partial y} \psi\, / \psi$,
where 
$\psi(x, y; a, \lambda) = a(1+x) + \left(e^{\lambda(x-1)} + e^{-\lambda(x-1)}\right)\cos(\lambda y)$
and $(a, \lambda)$ are parameters,
and its restriction to the boundary as Dirichlet boundary conditions.
The PDE is discretized using centered finite differences with with $482$ uniformly spaced grid points in the $x$-direction and $26$ uniformly spaced grid points in the $y$-direction. 
For ROM training, we collected $6400$ FOM snapshots corresponding to varying $(a, \lambda)\in [1, 10^4]\times [5, 25]$ (see Fig. \ref{fig:uv_states}) in a uniform $80\times 80$ grid. 
We use ROMs to predict the out-of-sample case $(a, \lambda) = (7692.5384, 21.9230)$.

First we use DD problem with $4$ uniformly sized subdomains in a $2\times 2$ configuration and vary the ROM sizes $n_i^\Omega$ and $n_i^\Gamma.$
Table \ref{tbl:error_speedup} shows that NM-ROM has an order of magnitude lower error than LS-ROM with and without HR when comparing ROMs of the same size.
In the non-HR case, LS-ROM only achieves order $10^{-3}$ error for a ROM with $96$ total DoF (error = $2.66 \times 10^{-3}$), while NM-ROM can achieve a similar error with only $36$ DoF (error = $2.42 \times 10^{-3}$) and a higher speedup (speedup = $26.2$) compared to LS-ROM with similar accuracy (speedup = $18.3$). 
LS-ROM achieves a much higher speedup in the HR cases while retaining similar errors from the non-HR cases. 
NM-ROM also retains high accuracy after HR, and gains an extra $15$-$20$ times speedup after applying HR.
\begin{table}[ht]
    \centering
    \begin{tabular}{|c|c|c|c|c|c|c|c|c|}
        \hline
         & $n_i^\Omega$ & $n_i^\Gamma$ & DoF & Error & Speedup & Error (HR) & Speedup (HR)\\
        \hline
        \multirow{4}{*}{LS-ROM} 
        &  $6$ & $3$ & $36$ & $2.06 \times 10^{-2} $ & $48.7$ & $1.78 \times 10^{-2} $ & $340.0$\\
        &  $8$ & $4$ & $48$ & $1.98 \times 10^{-2} $ & $30.0$ & $1.44 \times 10^{-2} $ & $347.6$\\
        & $10$ & $5$ & $60$ & $1.50 \times 10^{-2} $ & $16.3$ & $1.16 \times 10^{-2} $ & $329.6$\\
        & $16$ & $8$ & $96$ & $2.66 \times 10^{-3} $ & $18.3$ & $3.23 \times 10^{-3} $ & $280.4$\\
        \hline
        \multirow{4}{*}{NM-ROM} 
        &  $6$ & $3$ & $36$ & $2.42 \times 10^{-3} $ & $26.2$ & $2.60 \times 10^{-3} $ & $44.7$ \\
        &  $8$ & $4$ & $48$ & $1.28 \times 10^{-3} $ & $21.7$ & $1.64 \times 10^{-3} $ & $43.9$ \\
        & $10$ & $5$ & $60$ & $1.09 \times 10^{-3} $ & $15.0$ & $1.19 \times 10^{-3} $ & $43.6$ \\
        & $16$ & $8$ & $96$ & $7.87 \times 10^{-4} $ & $13.9$ & $9.80 \times 10^{-4} $ & $37.5$ \\
        \hline
      \end{tabular}
    \caption{Relative error and speedup for LS-ROM and NM-ROM with and without HR for varying ROM size. We use $N_i^B=100$ HR nodes per subdomain in the HR case.}
    \label{tbl:error_speedup}
\end{table}

Next we examine the per-subdomain reduction in the required number of autoencoder parameters for different subdomain configurations compared to the monolothic single-domain NM-ROM.
We use the notation $2\times 1$ subdomains to indicate $2$ subdomains in the $x$-direction and $1$ subdomain in the $y$-direction.
As expected, from Table \ref{tbl:NN_params}, we see that the maximum number of NN parameters per subdomain decreases significantly as more subdomains are used. 
Furthermore, the total number of NN parameters in the DD cases also decreases relative to the single-domain case. 
We also note that the error increases as more subdomains are used.
We kept the ROM size $(n_i^\Omega, n_i^\Gamma)=(6, 3)$ constant for each subdomain configuration to isolate the effect of DD on the number of NN parameters, but this may cause overfitting in the $16$ subdomain case. 
More careful hyper-parameter tuning is necessary to mitigate increases in error as the number of subdomains is increased. 

\begin{table}[ht]
    \centering
    \begin{tabular}{|c|c|c|c|c|}
        \hline
         Subdomains & Max \# subdomain params.\ & Reduction &  Total \# params. & Error\\
        \hline
          $1 \times 1$  & $2.995 \times 10^{6}$ &  $0.0$ \% & $2.995 \times 10^{6}$ & $1.08 \times 10^{-3}$\\
          \hline
          $2 \times 1$  & $1.147 \times 10^{6}$ & $61.7$ \% & $2.307 \times 10^{6}$ & $1.27 \times 10^{-3}$\\
          \hline
          $2 \times 2$  & $5.257 \times 10^{5}$ & $82.4$ \% & $2.384 \times 10^{6}$ & $2.42 \times 10^{-3}$\\
          \hline
          $4 \times 2$  & $2.617 \times 10^{5}$ & $91.3$ \% & $2.391 \times 10^{6}$ & $4.26 \times 10^{-3}$\\
          \hline
          $8 \times 2$  & $1.297 \times 10^{5}$ & $95.7$ \% & $2.406 \times 10^{6}$ & $4.58 \times 10^{-2}$\\
        \hline
      \end{tabular}
    \caption{Max number of NN parameters per subdomain, the per-subdomain reduction in number of NN parameters, the total number of parameters, and the corresponding error for different subdomain configurations. For the single-domain case, an NM-ROM of dimension $n=9$ is used. For the DD cases, $(n_i^\Omega, n_i^\Gamma)=(6, 3)$, resulting in $9$ DoF per subdomain. HR was not used to evaluate the NM-ROMs in these examples.}
    \label{tbl:NN_params}
\end{table}

\section{Conclusion}\label{sec:conclusion}
We extended the DD framework of \cite{CHoang_YChoi_KCarlberg_2021a}
and compute ROMs using NM-ROM with HR as presented in \cite{YKim_YChoi_DWidemann_TZohdi_2022a}.
Our experiments on the 2D Burgers' equation show that NM-ROM achieves an order of magnitude lower relative error than LS-ROM in nearly all cases tested. 
While LS-ROM with HR achieves much higher speedup than NM-ROM with HR, NM-ROM is still the clear winner in terms of ROM accuracy for a given ROM size. 
Moreover, HR allows NM-ROM to gain an extra $15$-$20$ time speedup compared to the non-HR cases. 
While the speedup is not as drastic as for LS-ROM, these speedup gains for NM-ROM are the highest that have been achieved for NM-ROM to our knowledge. 
We also showed that using the DD approach significantly decreases the number of required NN parameters per subdomain compared to the monolithic single-domain NM-ROM. 
In future work, we plan to apply DD NM-ROM to more challenging problems, including those with slowly decaying Kolmogorov $n$-width and to time-dependent problems. 
Other directions for future research include a greedy sampling strategy 
when choosing which FOM snapshots to compute for NM-ROM training and applying the DD NM-ROM framework to decomposable or component-based systems. 

\begin{ack}
This work was performed at Lawrence Livermore National Laboratory. 
A.\ N.\ Diaz was supported for this work by a Defense Science and Technology Internship (DSTI) at Lawrence Livermore National Laboratory and a 2021 National Defense Science and Engineering Graduate Fellowship.
Y.\ Choi was supported for this work by the U.S. Department of Energy, Office of Science, Office of Advanced Scientific Computing Research, as part of the CHaRMNET Mathematical Multifaceted Integrated Capability Center (MMICC) program, under Award Number DE-SC0023164 and partially by LDRD (21-SI-006).
M.\ Heinkenschloss was supported by AFOSR Grant FA9550-22-1-0004 at Rice University.
Lawrence Livermore National Laboratory is operated by Lawrence Livermore National Security, LLC, for the U.S. Department of Energy, National Nuclear Security Administration under Contract DE-AC52-07NA27344. IM review: LLNL-CONF-854737.
\end{ack}


\bibliography{references, yc_references}

\begin{thebibliography}{44}
\providecommand{\natexlab}[1]{#1}
\providecommand{\url}[1]{\texttt{#1}}
\expandafter\ifx\csname urlstyle\endcsname\relax
  \providecommand{\doi}[1]{doi: #1}\else
  \providecommand{\doi}{doi: \begingroup \urlstyle{rm}\Url}\fi

\bibitem[Antonietti et~al.(2016)Antonietti, Pacciarini, and Quarteroni]{PFAntonietti_PPacciarini_AQuarteroni_2016a}
P.~F. Antonietti, P.~Pacciarini, and A.~Quarteroni.
\newblock A discontinuous {G}alerkin reduced basis element method for elliptic problems.
\newblock \emph{ESAIM Math. Model. Numer. Anal.}, 50\penalty0 (2):\penalty0 337--360, 2016.
\newblock \doi{10.1051/m2an/2015045}.
\newblock URL \url{https://doi.org/10.1051/m2an/2015045}.

\bibitem[Antoulas(2005)]{ACAntoulas_2005a}
A.~C. Antoulas.
\newblock \emph{Approximation of Large-Scale Dynamical Systems}, volume~6 of \emph{Advances in Design and Control}.
\newblock Society for Industrial and Applied Mathematics (SIAM), Philadelphia, PA, 2005.
\newblock \doi{10.1137/1.9780898718713}.
\newblock URL \url{https://doi.org/10.1137/1.9780898718713}.

\bibitem[Antoulas et~al.(2016)Antoulas, Gosea, and Ionita]{ACAntoulas_IVGosea_ACIonita_2016a}
A.~C. Antoulas, I.~V. Gosea, and A.~C. Ionita.
\newblock Model reduction of bilinear systems in the {L}oewner framework.
\newblock \emph{SIAM J. Sci. Comput.}, 38\penalty0 (5):\penalty0 B889--B916, 2016.
\newblock URL \url{https://doi.org/10.1137/15M1041432}.

\bibitem[Antoulas et~al.(2020)Antoulas, Beattie, and Gugercin]{ACAntoulas_CABeattie_SGugercin_2020a}
A.~C. Antoulas, C.~A. Beattie, and S.~Gugercin.
\newblock \emph{Interpolatory Model Reduction}, volume~21 of \emph{Computational Science \& Engineering}.
\newblock Society for Industrial and Applied Mathematics (SIAM), Philadelphia, PA, 2020.
\newblock \doi{10.1137/1.9781611976083}.
\newblock URL \url{https://doi.org/10.1137/1.9781611976083}.

\bibitem[Barnett et~al.(2022)Barnett, Tezaur, and Mota]{JBarnett_ITezaur_AMota_2022a}
J.~Barnett, I.~Tezaur, and A.~Mota.
\newblock The {S}chwarz alternating method for the seamless coupling of nonlinear reduced order models and full order models.
\newblock \emph{arXiv:2210.12551}, 2022.
\newblock \doi{10.48550/ARXIV.2210.12551}.
\newblock URL \url{https://doi.org/10.48550/ARXIV.2210.12551}.

\bibitem[Barnett et~al.(2023)Barnett, Farhat, and Maday]{JBarnett_CFarhat_YMaday_2023a}
J.~Barnett, C.~Farhat, and Y.~Maday.
\newblock Neural-network-augmented projection-based model order reduction for mitigating the {K}olmogorov barrier to reducibility.
\newblock \emph{J. Comput. Phys.}, 492:\penalty0 Paper No. 112420, 20, 2023.
\newblock \doi{10.1016/j.jcp.2023.112420}.
\newblock URL \url{https://doi.org/10.1016/j.jcp.2023.112420}.

\bibitem[Benner and Breiten(2015)]{PBenner_TBreiten_2015a}
P.~Benner and T.~Breiten.
\newblock Two-sided projection methods for nonlinear model order reduction.
\newblock \emph{SIAM J. Sci. Comput.}, 37\penalty0 (2):\penalty0 B239--B260, 2015.
\newblock \doi{10.1137/14097255X}.
\newblock URL \url{http://dx.doi.org/10.1137/14097255X}.

\bibitem[Benner and Breiten(2017)]{PBenner_TBreiten_2017a}
P.~Benner and T.~Breiten.
\newblock Chapter 6: Model order reduction based on system balancing.
\newblock In P.~Benner, A.~Cohen, M.~Ohlberger, and K.~Willcox, editors, \emph{Model Reduction and Approximation: {T}heory and Algorithms}, Computational Science and Engineering, pages 261--295, Philadelphia, 2017. SIAM.
\newblock \doi{10.1137/1.9781611974829.ch6}.
\newblock URL \url{https://doi.org/10.1137/1.9781611974829.ch6}.

\bibitem[Buffoni et~al.(2009)Buffoni, Telib, and Iollo]{MBuffoni_HTelib_AIollo_2009a}
M.~Buffoni, H.~Telib, and A.~Iollo.
\newblock Iterative methods for model reduction by domain decomposition.
\newblock \emph{Comput. \& Fluids}, 38\penalty0 (6):\penalty0 1160--1167, 2009.
\newblock \doi{10.1016/j.compfluid.2008.11.008}.
\newblock URL \url{https://doi.org/10.1016/j.compfluid.2008.11.008}.

\bibitem[Carlberg et~al.(2018)Carlberg, Choi, and Sargsyan]{carlberg2018conservative}
K.~Carlberg, Y.~Choi, and S.~Sargsyan.
\newblock Conservative model reduction for finite-volume models.
\newblock \emph{Journal of Computational Physics}, 371:\penalty0 280--314, 2018.
\newblock \doi{10.1016/j.jcp.2018.05.019}.
\newblock URL \url{https://doi.org/10.1016/j.jcp.2018.05.019}.

\bibitem[Carlberg et~al.(2013)Carlberg, Farhat, Cortial, and Amsallem]{KTCarlberg_CFarhat_JCortial_DAmsallem_2013a}
K.~T. Carlberg, C.~Farhat, J.~Cortial, and D.~Amsallem.
\newblock The {GNAT} method for nonlinear model reduction: {E}ffective implementation and application to computational fluid dynamics and turbulent flows.
\newblock \emph{Journal of Computational Physics}, 242:\penalty0 623 -- 647, 2013.
\newblock \doi{10.1016/j.jcp.2013.02.028}.
\newblock URL \url{http://dx.doi.org/10.1016/j.jcp.2013.02.028}.

\bibitem[Cheung et~al.(2023)Cheung, Choi, Copeland, and Huynh]{cheung2023local}
S.~W. Cheung, Y.~Choi, D.~M. Copeland, and K.~Huynh.
\newblock Local lagrangian reduced-order modeling for the rayleigh-taylor instability by solution manifold decomposition.
\newblock \emph{Journal of Computational Physics}, 472:\penalty0 111655, 2023.
\newblock \doi{10.1016/j.jcp.2022.111655}.
\newblock URL \url{https://doi.org/10.1016/j.jcp.2022.111655}.

\bibitem[Choi and Carlberg(2019)]{choi2019space}
Y.~Choi and K.~Carlberg.
\newblock Space--time least-squares petrov--galerkin projection for nonlinear model reduction.
\newblock \emph{SIAM Journal on Scientific Computing}, 41\penalty0 (1):\penalty0 A26--A58, 2019.
\newblock \doi{10.1137/17M1120531}.
\newblock URL \url{https://doi.org/10.1137/17M1120531}.

\bibitem[Choi et~al.(2021)Choi, Brown, Arrighi, Anderson, and Huynh]{choi2021space}
Y.~Choi, P.~Brown, W.~Arrighi, R.~Anderson, and K.~Huynh.
\newblock Space--time reduced order model for large-scale linear dynamical systems with application to boltzmann transport problems.
\newblock \emph{Journal of Computational Physics}, 424:\penalty0 109845, 2021.
\newblock \doi{10.1016/j.jcp.2020.109845}.
\newblock URL \url{https://doi.org/10.1016/j.jcp.2020.109845}.

\bibitem[Copeland et~al.(2022)Copeland, Cheung, Huynh, and Choi]{copeland2022reduced}
D.~M. Copeland, S.~W. Cheung, K.~Huynh, and Y.~Choi.
\newblock Reduced order models for lagrangian hydrodynamics.
\newblock \emph{Computer Methods in Applied Mechanics and Engineering}, 388:\penalty0 114259, 2022.
\newblock \doi{10.1016/j.cma.2021.114259}.
\newblock URL \url{https://doi.org/10.1016/j.cma.2021.114259}.

\bibitem[Diaz et~al.(2023)Diaz, Choi, and Heinkenschloss]{ANDiaz_YChoi_MHeinkenschloss_2023a}
A.~N. Diaz, Y.~Choi, and M.~Heinkenschloss.
\newblock A fast and accurate domain-decomposition nonlinear manifold reduced order model.
\newblock \emph{arXiv:2305.15163v1}, 2023.
\newblock \doi{10.48550/arXiv.2305.15163}.
\newblock URL \url{https://doi.org/10.48550/arXiv.2305.15163}.

\bibitem[Eftang and Patera(2013)]{JLEftang_ATPatera_2013a}
J.~L. Eftang and A.~T. Patera.
\newblock Port reduction in parametrized component static condensation: approximation and {\it a posteriori} error estimation.
\newblock \emph{Internat. J. Numer. Methods Engrg.}, 96\penalty0 (5):\penalty0 269--302, 2013.
\newblock \doi{10.1002/nme.4543}.
\newblock URL \url{https://doi.org/10.1002/nme.4543}.

\bibitem[Eftang et~al.(2012)Eftang, Huynh, Knezevic, Ronquist, and Patera]{JLEftang_DBPHuynh_DJKnezevic_EMRonquist_ATPatera_2012a}
J.~L. Eftang, D.~B.~P. Huynh, D.~J. Knezevic, E.~M. Ronquist, and A.~T. Patera.
\newblock Adaptive port reduction in static condensation.
\newblock \emph{IFAC Proceedings Volumes}, 45\penalty0 (2):\penalty0 695--699, 2012.
\newblock \doi{10.3182/20120215-3-AT-3016.00123}.
\newblock URL \url{https://doi.org/10.3182/20120215-3-AT-3016.00123}.
\newblock 7th Vienna International Conference on Mathematical Modelling.

\bibitem[Gosea and Antoulas(2018)]{IVGosea_ACAntoulas_2018a}
I.~V. Gosea and A.~C. Antoulas.
\newblock Data-driven model order reduction of quadratic-bilinear systems.
\newblock \emph{Numer. Linear Algebra Appl.}, 25\penalty0 (6):\penalty0 e2200, 2018.
\newblock \doi{10.1002/nla.2200}.
\newblock URL \url{http://dx.doi.org/10.1002/nla.2200}.

\bibitem[Gu(2011)]{CGu_2011a}
C.~Gu.
\newblock {QLMOR}: A projection-based nonlinear model order reduction approach using quadratic-linear representation of nonlinear systems.
\newblock \emph{Computer-Aided Design of Integrated Circuits and Systems, IEEE Transactions on}, 30\penalty0 (9):\penalty0 1307--1320, sept. 2011.
\newblock \doi{10.1109/TCAD.2011.2142184}.
\newblock URL \url{https://doi.org/10.1109/TCAD.2011.2142184}.

\bibitem[Gubisch and Volkwein(2017)]{MGubisch_SVolkwein_2017a}
M.~Gubisch and S.~Volkwein.
\newblock Chapter 1: {P}roper {O}rthogonal {D}ecomposition for linear-quadratic optimal control.
\newblock In P.~Benner, A.~Cohen, M.~Ohlberger, and K.~Willcox, editors, \emph{Model Reduction and Approximation: {T}heory and Algorithms}, Computational Science and Engineering, pages 3--64, Philadelphia, 2017. SIAM.
\newblock \doi{10.1137/1.9781611974829.ch1}.
\newblock URL \url{https://doi.org/10.1137/1.9781611974829.ch1}.

\bibitem[Haasdonk(2017)]{BHaasdonk_2017a}
B.~Haasdonk.
\newblock Chapter 2: Reduced basis methods for parametrized {PDE}s - a tutorial introduction for stationary and instationary problems.
\newblock In P.~Benner, A.~Cohen, M.~Ohlberger, and K.~Willcox, editors, \emph{Model Reduction and Approximation: {T}heory and Algorithms}, Computational Science and Engineering, pages 65--136. SIAM, Philadelphia, 2017.
\newblock \doi{10.1137/1.9781611974829.ch2}.
\newblock URL \url{https://doi.org/10.1137/1.9781611974829.ch2}.

\bibitem[Hartman and Mestha(2017)]{DHartman_LKMestha_2017a}
D.~Hartman and L.~K. Mestha.
\newblock A deep learning framework for model reduction of dynamical systems.
\newblock In \emph{2017 IEEE Conference on Control Technology and Applications (CCTA)}, pages 1917--1922, 2017.
\newblock \doi{10.1109/CCTA.2017.8062736}.
\newblock URL \url{https://doi.org/10.1109/CCTA.2017.8062736}.

\bibitem[Hinze and Volkwein(2005)]{MHinze_SVolkwein_2005a}
M.~Hinze and S.~Volkwein.
\newblock Proper orthogonal decomposition surrogate models for nonlinear dynamical systems: Error estimates and suboptimal control.
\newblock In P.~Benner, V.~Mehrmann, and D.~C. Sorensen, editors, \emph{Dimension Reduction of Large-Scale Systems}, Lecture Notes in Computational Science and Engineering, Vol.~45, pages 261--306, Heidelberg, 2005. Springer-Verlag.
\newblock \doi{10.1007/3-540-27909-1_10}.
\newblock URL \url{http://doi.org/10.1007/3-540-27909-1_10}.

\bibitem[Hoang et~al.(2021)Hoang, Choi, and Carlberg]{CHoang_YChoi_KCarlberg_2021a}
C.~Hoang, Y.~Choi, and K.~Carlberg.
\newblock Domain-decomposition least-squares {P}etrov-{G}alerkin ({DD}-{LSPG}) nonlinear model reduction.
\newblock \emph{Comput. Methods Appl. Mech. Engrg.}, 384:\penalty0 Paper No. 113997, 41, 2021.
\newblock \doi{10.1016/j.cma.2021.113997}.
\newblock URL \url{https://doi.org/10.1016/j.cma.2021.113997}.

\bibitem[Huynh et~al.(2013)Huynh, Knezevic, and Patera]{DBPHuynh_DJKnezevic_ATPatera_2013a}
D.~B.~P. Huynh, D.~J. Knezevic, and A.~T. Patera.
\newblock A static condensation reduced basis element method: approximation and {\it a posteriori} error estimation.
\newblock \emph{ESAIM Math. Model. Numer. Anal.}, 47\penalty0 (1):\penalty0 213--251, 2013.
\newblock \doi{10.1051/m2an/2012022}.
\newblock URL \url{https://doi.org/10.1051/m2an/2012022}.

\bibitem[Iapichino et~al.(2012)Iapichino, Quarteroni, and Rozza]{LIapichino_AQuarteroni_GRozza_2012a}
L.~Iapichino, A.~Quarteroni, and G.~Rozza.
\newblock A reduced basis hybrid method for the coupling of parametrized domains represented by fluidic networks.
\newblock \emph{Comput. Methods Appl. Mech. Engrg.}, 221/222:\penalty0 63--82, 2012.
\newblock \doi{10.1016/j.cma.2012.02.005}.
\newblock URL \url{https://doi.org/10.1016/j.cma.2012.02.005}.

\bibitem[Iollo et~al.(2023)Iollo, Sambataro, and Taddei]{AIollo_GSambataro_TTaddei_2022a}
A.~Iollo, G.~Sambataro, and T.~Taddei.
\newblock A one-shot overlapping {S}chwarz method for component-based model reduction: application to nonlinear elasticity.
\newblock \emph{Comput. Methods Appl. Mech. Engrg.}, 404:\penalty0 Paper No. 115786, 32, 2023.
\newblock \doi{10.1016/j.cma.2022.115786}.
\newblock URL \url{https://doi.org/10.1016/j.cma.2022.115786}.

\bibitem[Kashima(2016)]{KKashima_2016a}
K.~Kashima.
\newblock Nonlinear model reduction by deep autoencoder of noise response data.
\newblock In \emph{2016 IEEE 55th Conference on Decision and Control (CDC)}, pages 5750--5755, 2016.
\newblock \doi{10.1109/CDC.2016.7799153}.
\newblock URL \url{https://doi.org/10.1109/CDC.2016.7799153}.

\bibitem[Kim et~al.(2020)Kim, Choi, Widemann, and Zohdi]{kim2020efficient}
Y.~Kim, Y.~Choi, D.~Widemann, and T.~Zohdi.
\newblock Efficient nonlinear manifold reduced order model.
\newblock \emph{arXiv preprint arXiv:2011.07727}, 2020.
\newblock \doi{10.48550/arXiv.2011.07727}.
\newblock URL \url{https://doi.org/10.48550/arXiv.2011.07727}.

\bibitem[Kim et~al.(2021)Kim, Wang, and Choi]{kim2021efficient}
Y.~Kim, K.~Wang, and Y.~Choi.
\newblock Efficient space--time reduced order model for linear dynamical systems in python using less than 120 lines of code.
\newblock \emph{Mathematics}, 9\penalty0 (14):\penalty0 1690, 2021.
\newblock \doi{10.3390/math9141690}.
\newblock URL \url{https://doi.org/10.3390/math9141690}.

\bibitem[Kim et~al.(2022)Kim, Choi, Widemann, and Zohdi]{YKim_YChoi_DWidemann_TZohdi_2022a}
Y.~Kim, Y.~Choi, D.~Widemann, and T.~Zohdi.
\newblock A fast and accurate physics-informed neural network reduced order model with shallow masked autoencoder.
\newblock \emph{J. Comput. Phys.}, 451:\penalty0 Paper No. 110841, 29, 2022.
\newblock \doi{10.1016/j.jcp.2021.110841}.
\newblock URL \url{https://doi.org/10.1016/j.jcp.2021.110841}.

\bibitem[Lee and Carlberg(2020)]{KLee_KTCarlberg_2020a}
K.~Lee and K.~T. Carlberg.
\newblock Model reduction of dynamical systems on nonlinear manifolds using deep convolutional autoencoders.
\newblock \emph{J. Comput. Phys.}, 404:\penalty0 108973, 32, 2020.
\newblock \doi{10.1016/j.jcp.2019.108973}.
\newblock URL \url{https://doi.org/10.1016/j.jcp.2019.108973}.

\bibitem[Li et~al.(2020{\natexlab{a}})Li, Tang, Wu, and Liao]{KLi_KTang_TWu_QLiao_2020a}
K.~Li, K.~Tang, T.~Wu, and Q.~Liao.
\newblock {D3M}: A deep domain decomposition method for partial differential equations.
\newblock \emph{IEEE Access}, 8:\penalty0 5283--5294, 2020{\natexlab{a}}.
\newblock \doi{10.1109/ACCESS.2019.2957200}.
\newblock URL \url{https://doi.org/10.1109/ACCESS.2019.2957200}.

\bibitem[Li et~al.(2023)Li, Xia, Liu, and Liao]{SLi_YXia_YLiu_QLiao_2023a}
S.~Li, Y.~Xia, Y.~Liu, and Q.~Liao.
\newblock A deep domain decomposition method based on fourier features.
\newblock \emph{Journal of Computational and Applied Mathematics}, 423:\penalty0 114963, 2023.
\newblock \doi{10.1016/j.cam.2022.114963}.
\newblock URL \url{https://doi.org/10.1016/j.cam.2022.114963}.

\bibitem[Li et~al.(2020{\natexlab{b}})Li, Xiang, and Xu]{WLi_XXiang_YXu_2020a}
W.~Li, X.~Xiang, and Y.~Xu.
\newblock Deep domain decomposition method: {E}lliptic problems.
\newblock In J.~Lu and R.~Ward, editors, \emph{Proceedings of The First Mathematical and Scientific Machine Learning Conference}, volume 107 of \emph{Proceedings of Machine Learning Research}, pages 269--286. PMLR, 20--24 Jul 2020{\natexlab{b}}.
\newblock URL \url{https://proceedings.mlr.press/v107/li20a.html}.

\bibitem[Maday and R{\o}nquist(2002)]{YMaday_EMRonquist_2002a}
Y.~Maday and E.~M. R{\o}nquist.
\newblock A reduced-basis element method.
\newblock \emph{J. Sci. Comput.}, 17\penalty0 (1-4):\penalty0 447--459, 2002.
\newblock \doi{10.1023/A:1015197908587}.
\newblock URL \url{https://doi.org/10.1023/A:1015197908587}.

\bibitem[Maday and R{\o}nquist(2004)]{YMaday_EMRonquist_2004a}
Y.~Maday and E.~M. R{\o}nquist.
\newblock The reduced basis element method: application to a thermal fin problem.
\newblock \emph{SIAM J. Sci. Comput.}, 26\penalty0 (1):\penalty0 240--258, 2004.
\newblock \doi{10.1137/S1064827502419932}.
\newblock URL \url{https://doi.org/10.1137/S1064827502419932}.

\bibitem[Mayo and Antoulas(2007)]{AJMayo_ACAntoulas_2007a}
A.~J. Mayo and A.~C. Antoulas.
\newblock A framework for the solution of the generalized realization problem.
\newblock \emph{Linear Algebra Appl.}, 425\penalty0 (2-3):\penalty0 634--662, 2007.
\newblock \doi{10.1016/j.laa.2007.03.008}.
\newblock URL \url{https://doi.org/10.1016/j.laa.2007.03.008}.

\bibitem[Ohlberger and Rave(2016)]{MOhlberger_SRave_2016a}
M.~Ohlberger and S.~Rave.
\newblock Reduced basis methods: Success, limitations and future challenges.
\newblock \emph{Proceedings of the Conference Algoritmy}, pages 1--12, 2016.
\newblock URL \url{http://www.iam.fmph.uniba.sk/amuc/ojs/index.php/algoritmy/article/view/389}.

\bibitem[Prechelt(1998)]{LPrechelt_1998a}
L.~Prechelt.
\newblock Automatic early stopping using cross validation: quantifying the criteria.
\newblock \emph{Neural networks}, 11\penalty0 (4):\penalty0 761--767, 1998.
\newblock \doi{10.1016/S0893-6080(98)00010-0}.
\newblock URL \url{https://doi.org/10.1016/S0893-6080(98)00010-0}.

\bibitem[Quarteroni et~al.(2016)Quarteroni, Manzoni, and Negri]{AQuarteroni_AManzoni_FNegri_2016a}
A.~Quarteroni, A.~Manzoni, and F.~Negri.
\newblock \emph{Reduced Basis Methods for Partial Differential Equations. {A}n Introduction}, volume~92 of \emph{Unitext}.
\newblock Springer, Cham, 2016.
\newblock \doi{10.1007/978-3-319-15431-2}.
\newblock URL \url{https://doi.org/10.1007/978-3-319-15431-2}.

\bibitem[Smetana and Taddei(2022)]{KSmetana_TTaddei_2022a}
K.~Smetana and T.~Taddei.
\newblock Localized model reduction for nonlinear elliptic partial differential equations: localized training, partition of unity, and adaptive enrichment.
\newblock \emph{arXiv:2202.09872v1}, 2022.
\newblock \doi{10.48550/ARXIV.2202.09872}.
\newblock URL \url{https://doi.org/10.48550/ARXIV.2202.09872}.

\bibitem[Sun et~al.(2022)Sun, Xu, and Yi]{QSun_XXu_HYi_2023a}
Q.~Sun, X.~Xu, and H.~Yi.
\newblock Domain decomposition learning methods for solving elliptic problems.
\newblock \emph{arXiv preprint arXiv:2207.10358}, 2022.
\newblock \doi{10.48550/arXiv.2207.10358}.
\newblock URL \url{https://doi.org/10.48550/arXiv.2207.10358}.

\end{thebibliography}


\end{document}